\documentclass[10pt]{article}
\usepackage{amssymb}
\usepackage{amsmath}
\usepackage[dvips]{graphics}
\usepackage{epsfig}
\usepackage{pgfplots}
\usepackage{pgfplotstable}
\usepackage{float}
\pgfplotsset{compat=1.17}

\newtheorem{theorem}{Theorem}

\newtheorem{proposition}[theorem]{Proposition}

\begin{document}
\date{}
\author{Aristides V. Doumas\footnote{
National Technical University of Athens, Department of Mathematics, Zografou Campus,
157 80 Athens, GREECE, and Archimedes/Athena Research Center, GREECE.
adou@math.ntua.gr, aris.doumas@hotmail.com}}
\title{Dixie cup problem in an interlacing process}
\maketitle

\centerline{\textit{Dedicated to the memory of professor Miltos Makris.}}
\centerline{\textit{Always a cadet, brother, and Best Friend.}}
\begin{abstract}
The “double Dixie cup problem” of D.J. Newman and L. Shepp is a well-known variant of the coupon collector’s problem, where the object of study is the number of coupons that a collector has to buy in order to complete $m$ sets of all $N$ existing different coupons.
In this paper we consider the case where the coupons distribution is a \textit{mixture} of two different distributions, where the coupons from the first distribution are far rarer than the ones coming from the second. We apply a Poissonization technique, as well as well known results and techniques from our previous work, to derive the asymptotics (leading term) of the expectation of the above random variable as $N\rightarrow \infty$ for large classes of distributions. As it turns out, both distributions contribute to this result. The leading asymptotics of the rising moments of the aforementioned random variable are also discussed. We conclude by generalizing the problem to the case where the family of coupons is a mixture of $j$ subfamilies.
\end{abstract}

\textbf{MSC 2022 Mathematics Classification.} 40-08; 40G99; 65A16.

\section{Introduction}
The “double Dixie cup problem” of D.J. Newman and L. Shepp \cite{NS} refers to a population whose members are of $N$ different types (coupons).\footnote{It is entertaining that the famous American girl group The Dixie Cups was active at the beginning of the 1960s, i.e., around the time of \cite{NS}. } For $j = 1, 2, . . . , N$, we denote by $p_{j}$ the probability that a randomly chosen coupon is of type $j$, where $p_{j}>0$ and $\sum_{j=1}^{N}p_{j}=1$. The members of the population are sampled independently with replacement and their types are recorded. Let $\alpha=\left\{a_{j}\right\}_{j=1}^{\infty}$ be a sequence of strictly positive numbers. Then, for each positive integer $N$
one can create a probability measure
$\pi _N =\{p_1,...,p_N\}$ on the set of types $\{1,...,N\}$ by taking
\begin{equation}
p_j = \frac{a_j}{A_N},\,\,\, \text{where}\,\,\,
A_N = \sum_{j=1}^N a_j.
\label{0}
\end{equation}
The number $T_{m}\left(N;\alpha\right)$ of trials needed until all $N$ types are detected (at least $m$ times) is, naturally, of interest. In this paper we consider the case where the sequence $\alpha$ is formed as the \textit{interlacing} of two subsequences and derive the asymptotics (leading term) of the expectation of the random variable $T_{m}\left(N;\alpha\right)$ as $N\rightarrow \infty$.\\
Applications of the classic problem arise in many areas of science, such as computer science (search algorithms), physics, biology, ecology, earth and planetary sciences, economics and finance, as well as linguistics, demography, and the social sciences; see, e.g., \cite{BH}. The variant of the problem studied here is often applied in the same contexts, since there are cases where half of the coupons follow a different distribution (and are in fact rarer) than the other half. In cybersecurity, the Dixie cup problem can model the process of discovering all vulnerabilities in a network. Each unique vulnerability corresponds to a coupon, and each security scan to an attempt to collect one. Since not all vulnerabilities are equally likely, it is natural to distinguish two categories: common vulnerabilities, which are easily detected by standard scanning tools, and rare vulnerabilities, which are harder to uncover and follow a separate probability distribution. This framework allows one to predict the number of scans required to achieve a complete collection (the case where $m=1$), i.e., to identify every vulnerability in the system. A similar application arises in epidemiology, where scientists track the emergence of virus strains within a population. Some strains are common, while others are rare or geographically localized, and the process of collecting samples to detect all strains can likewise be modeled by this problem.\\
We now return to the above-mentioned random variable $T_{m}\left(N;\alpha\right)$. Its average w.r.t. to the sequence $\alpha$ will be denoted by $E\left[T_{m}\left(N;\alpha\right)\right]$. 
It is well known (see, e.g., \cite{DPSETS}) that
\begin{equation}
E\left[T_{m}\left(N;\alpha\right)\right] =
 \int_{0}^{\infty}\left[1-\prod_{j=1}^N
 \bigg(1-S_m({p_{j}t})e^{-{p_{j}}t}\bigg)\right]dt,\label{bsc}
\end{equation}
where $S_m (y)$ is the $m-th$ partial sum of $e^{y}$, namely
\begin{equation*}
  S_m(y)=\sum_{i=0}^{m-1}\frac{y^{i}}{i!}.
\end{equation*}
From now on we assume that $N$ is even and consider the real sequences $\beta:=\{b_j\}_{j=1}^{\infty }$ and $\delta:=\{d_j\}_{j=1}^{\infty }$ such that
\begin{equation}
a_{2j-1} = b_j
\qquad \text{and} \qquad
a_{2j} =  d_j,
\qquad\qquad
j = 1, 2, \dots.
\label{sb1}
\end{equation}
For convenience we set
\begin{equation}
N := 2M.\label{nm}
\end{equation}
Hence,
\begin{align}
A_N =B_M +D_M, \label{01}
\end{align}
where
\begin{align}
B_M:=\sum_{j=1}^M b_j\,\,\,\,\,\text{and}\,\,\,\,D_M:=\sum_{j=1}^M d_j.\label{011}
\end{align}
In words, the set of $N$ coupons consists of two different sets, each one having $M$ coupons with the corresponding coupon probabilities
\begin{equation}
  q_{j}=\frac{b_{j}}{B_M}\,\,\,\text{and}\,\,\,\, \tilde{q_{j}}=\frac{d_{j}}{D_M}.\label{stdccp}
\end{equation}
Notice that, when $M\rightarrow \infty$, then, of course, $N\rightarrow \infty$; thanks to (\ref{nm}). In this paper we consider \textit{growing} sequences $\beta:=\{b_j\}_{j=1}^{\infty }$, such that
\begin{equation}
  b_{j}\rightarrow \infty\label{deca}
\end{equation}
and \textit{decaying} sequences $\delta:=\{d_j\}_{j=1}^{\infty }$ of the form
\begin{equation}
d_{j}=\frac{1}{f(j)},\label{cond2}
\end{equation}
where
\begin{equation}
f(x)>0\qquad \text{and}\qquad f^{\prime}(x)>0,\qquad x>0. \label{cond2a}
\end{equation}
Furthermore we assume that $f(x)$ possesses three derivatives and satisfies
the following conditions as $x\rightarrow \infty$:
\begin{align}
\text{(i)}\,\,f(x)\rightarrow \infty,\,\,\text{(ii)}\,\, \frac{f^{\prime }(x)}{f(x)}\rightarrow 0,\,\,\text{(iii)} \frac{f^{\prime \prime}(x)/f^{\prime }(x)}{f'(x)/f(x)} = O\left(1\right),\,\,\text{(iv) } \frac{f^{\prime \prime\prime}(x)\;f(x)^{2}}{ f^{\prime }(x)^{3}} = O\left(1\right).\label{cond3}
\end{align}
In words $f(\cdot)$ belongs to the class of positive and strictly increasing $C^{3}\left(0,\infty\right)$ functions, which grow to $\infty$ (as $x\rightarrow \infty$) slower than exponentials, but faster than powers of logarithms. Examples of the growing sequences $\beta$ for are: $j^{p},\,p>0,\,\,\,e^{pj},\,p>0,\,\,\,j!$, and for the sequences $\delta$ are: $j^{-p}\left(\ln j\right)^{-q},\,p>0,\,q\in \mathbb{R},\,\,e^{-j^{r}}$, where, $0<r<1$. This work \textit{also covers} the case where one of the two sequences $\beta$ or $\delta$ is the uniform distribution. Recently and just for the case where $m=1$, the scenario where the sequence $\alpha$ is the mixing of two subsequences one of which is constant (this corresponds to a uniform subcollection of coupons), while the other obeys the well-known \textit{Zipf} law has been considered; we will come back to this case later. Hence, we may say that the coupons produced by the sequence $\delta$ are \textit{far rarer} than the ones from the sequence $\beta$. We conclude this section by mentioning that $T_{m}\left(M;\beta\right)$ denotes the number of trials needed to obtain $m$ complete sets of the $M$ coupons generated by the sequence $\beta$ via relation (\ref{stdccp}). In other words, in the Dixie cup problem concerning $m$ complete sets of $M$ different coupons, each coupon has probability  $q_{j}=b_{j} /B_M$. (Similarly, one defines  $T_{m}\left(M;\delta\right)$). 

\section{Asymptotics of $T_{m}(N; \alpha)$}
Motivated by the seminal work of Holst (see, \cite{HOLST}) we consider a homogeneous Poisson process $\left\{N_{t}\right\}_{t\ge 0}$ with intensity one. Each event of this process is a collected coupon, and is classified as type $j$ with probability $p_{j}$. Set 
\begin{equation*}
 N_t:= X_t +Y_t,
\end{equation*}
where $\left\{X_{t}\right\}_{t\ge 0}$ and $\left\{Y_{t}\right\}_{t\ge 0}$ are two independent Poisson processes, such that each Poisson event of them is a collected coupon coming from the sequences $\beta$ and $\delta$ respectively, thanks to (\ref{01}) and (\ref{011}). (This  Poissonization technique is also presented in detail in \cite{SR}). Recall that from our assumptions the coupons coming from the sequence $\delta$ are sufficiently rarer than the ones coming from the sequence $\beta$. Hence, as $N\rightarrow \infty$ with probability 1 the last coupon collected will be from the sequence $\delta$.
Let $W_{j}$ denote the random variable representing the number of arrivals of the process $\left\{X_{t}\right\}_{t\ge 0}$ between two consecutive arrivals of $\left\{Y_{t}\right\}_{t\ge 0}$. Then, $W_{j}$'s are i.i.d random variables. Moreover, it is an easy exercise for one to see that they are geometrically distributed with parameter $\frac{D_M}{A_N}$, namely:
\begin{equation}
  W_{j}\sim \text{Geo}\left(\frac{D_M}{A_N}\right), \,\,\,\,\,j=1,2,\dots,\label{geo}
\end{equation}
(where we have used (\ref{01})). Hence, one has
\begin{equation*}
  T_{m}\left(N;\alpha\right)=\sum_{j=1}^{T_{m}\left(M;\delta\right)}W_{j},
  \,\,\,\,\,\text{a.s.},\qquad N\rightarrow \infty.
\end{equation*}
It follows that
\begin{align}
  E\left[T_{m}\left(N;\alpha\right)\right]=&
  E\left[T_{m}\left(M;\delta\right)\right]E\left[W_{j}\right] \nonumber\\
  =&E\left[T_{m}\left(M;\delta\right)\right] \left(\frac{D_M +B_M}{D_M}\right) \nonumber\\
  =&E\left[T_{m}\left(M;\delta\right)\right] \left(\frac{A_N}{D_M}\right),\label{02}
\end{align}
where we have used Wald's lemma in the first equation above (since the finite mean random variable $T_{m}\left(M;\delta\right)$ is a stopping time with respect to the finite mean sequence $\left\{W_{j}\right\}_{j \in \mathbb{N} }$), the distribution of $W_{j}$'s in the second, and relation (\ref{01}) in the last one. Since,
\begin{equation*}
E\left[T_{m}\left(M;\delta\right)\right] =
 \int_{0}^{\infty}\left[1-\prod_{j=1}^M
 \bigg(1-S_m(\tilde{q}_{j}t)e^{-\tilde{q_{j}}\,\,t}\bigg)\right]dt,
\end{equation*}
and $\tilde{q}_{j}=\frac{d_{j}}{D_M}$ the change of variables $t=u D_M$ yields
\begin{equation}
E\left[T_{m}\left(M;\delta\right)\right] =
D_M \int_{0}^{\infty}\left[1-\prod_{j=1}^M
 \bigg(1-S_m(d_{j}u)e^{-d_{j}u}\bigg)\right]du.\label{1987}
\end{equation}
By invoking (\ref{1987}) in (\ref{02}) one has proved the following
\begin{theorem}
Let the sequence $\alpha = \{a_{j}\}_{j=1}^{\infty }$ is formed as the \textit{interlacing} of two subsequences $\beta = \{b_{j}\}_{j=1}^{\infty }$ and $\delta = \{d_{j}\}_{j=1}^{\infty }$, as given in relation (\ref{sb1}).  Then,  as $N=2M \rightarrow \infty$ we have
\begin{equation}
E\left[T_{m}\left(N;\alpha\right)\right]\sim \left(D_M +B_M\right)\int_{0}^{\infty}\left[1-\prod_{j=1}^M
 \bigg(1-S_m(d_{j}u)e^{-d_{j}u}\bigg)\right]du \label{03}
\end{equation}
where, 
\begin{align*}
B_M:=\sum_{j=1}^M b_j\,\,\,\,\,\text{and}\,\,\,\,D_M:=\sum_{j=1}^M d_j,
\end{align*}
provided that the coupons from the sequence $\delta$ are far rarer than the ones from the sequence $\beta$ (see, (\ref{deca})--(\ref{cond3})).
\end{theorem}
From here and in what follows, the asymptotic notation used in relation (\ref{03}) has the meaning that $a_n \sim b_n$, as $n\rightarrow \infty$ means that $\lim_{n\rightarrow \infty} a_n /b_n =1$. Relation (\ref{03}) is the \textit{key formula} to derive leading asymptotics for $E\left[T_{m}\left(N;\alpha\right)\right]$ as $N \rightarrow \infty$. It follows that our problem can be treated as \textit{two
separate problems}, namely gaining inside the order of magnitude of
\begin{equation*}
  A_N=D_M +B_M
\end{equation*}
and approximating the integral appearing in (\ref{03}), i.e.:
\begin{equation*}
L_{1}\left(M;\delta\right):=\int_{0}^{\infty}\left[1-\prod_{j=1}^M
 \bigg(1-S_m(d_{j}u)\bigg)e^{-d_{j}u}\right]du.
\end{equation*}
As we will see later, the number $m$ \textit{does not} appear in the leading term of the average as $M\rightarrow \infty$. Intitutevely, $m$ is invisible in the leading order, since the first few coupons are easy to collect, while the last few missing coupon types (the “hard part”) dominates the total time. \\
\textbf{Step I}. The good news is that asymptotics of the quantity $A_N$ (as $M \rightarrow \infty$, i.e., as $N \rightarrow \infty$) can be can be handled by existing powerful methods, such as the celebrated  Euler-Maclaurin summation formula (see, e.g., \cite{B-O}). In particular, the leading term in the asymptotic expansion of the sum $\sum_{j=1}^{M}l_j$ (for some sequence $\lambda:=\left\{l_{j}\right\}_{j=1}^{\infty}$) depends on the behaviour of the series $\sum_{j=1}^{\infty}l_j$. If
\begin{equation*}
C_\lambda :=\sum_{j=1}^{\infty}l_j<\infty,
\end{equation*}
then
\begin{equation*}
\sum_{j=1}^{M}l_j \sim C_\lambda ,\,\,\, \,\,\,M\rightarrow \infty,
\end{equation*}
while, in the case where
\begin{equation*}
C_\lambda=\infty,
\end{equation*}
we have
\begin{equation}
\,\,\,\,\,\,\,\,\,\,\,\,\,\,\,\,\,\,\,\sum_{j=1}^{M}l_j \sim \int_{1}^{M}l(x)dx,\,\,\,\,\, \,\,\,M\rightarrow \infty.\label{05}
\end{equation}
where  $l(x)$ is the real function associated with the sequence $\left\{l_{j}\right\}_{j=1}^{\infty}$. To mention a few examples we have (as $M\rightarrow \infty$):
\begin{align*}
\,\,\,\,\,\,\,\,\,\,\,\,\,\,\,\,\,\,\,
\,\,\,\,\,\,\,\,\,\,\,\,\,\,\,\,\,\,\,\,\,\,\,\,\,\,\,\,\,\sum_{j=1}^{M}j^{p}&\sim \frac{M^{p+1}}{p+1},\,\,\,\,p>0,\text{(positive power law)},  \end{align*}
and
\begin{equation*}
\,\,\,\,\,\,\,\,\,\,\,\,\,\,\,\,\,\,\,
\,\,\,\,\,\,\,\,\,\,\,\,\,\,\,\,\,\,\,\,\,\,\,\,\,\,\,\,\,\,\,\,\,
\,\,\,\,\,\,\,\,\,\,\,\,\,\,\,
\sum_{j=1}^{M}j^{-p} \sim
\left\{
	\begin{array}{ll}
		\frac{M^{1-p}}{1-p}  & \mbox{if }\,\, p\in (0,1), \\
        \ln M  & \mbox{if }\,\, p=1,\,\,\,\,\,\,\,\,\,\text{(generalized Zipf law)} \\		\zeta(p) & \mbox{if }\,\, p>1,
	\end{array}
\right.
\end{equation*}
where $\zeta(\cdot)$ is the Riemann’s Zeta function. Another interesting example is the one where $l_{j}=\left(\ln j\right)^{-p},\,\,p>0$. Here, $C_\lambda = \infty$. Hence, (\ref{05}) and integration by parts yields
\begin{equation*}
  \sum_{j=2}^{M}l_j \sim \frac{M}{\left(\ln M\right)^p},\,\,\,\,\,\,\,\,M\rightarrow \infty.
\end{equation*}
In particular, for $p = 1$ we get the so-called \textit{offset logarithmic integral} or \textit{Eulerian logarithmic integral}, namely,
\begin{equation*}
Li(M):=\int_{2}^{M}\frac{dt}{\ln t}
\end{equation*}
which is a very good approximation to the number of prime numbers less than $M$ (as $M\rightarrow \infty$).\\
Let us now consider the case where  $l_{j}=\frac{1}{j!},\, j=0,1,\cdots .$ Since, $\sum_{j=0}^{\infty}\frac{1}{j!}=e$ we have
\begin{equation*}
  \sum_{j=0}^{M}\frac{1}{j!}\sim e, \,\,\,\,\,\,\,\,M\rightarrow \infty.
\end{equation*}
\textbf{Remark 1}. Notice that there are cases where the Euler-Maclaurin summation formula is \textit{not} effective. For example, if $l_{j}=e^{jc_j}$, where $c_j$ is increasing and $\lim_{j\rightarrow \infty}c_j =\infty$. Then,
\begin{equation}
  \sum_{j=1}^{M}e^{jc_j} \sim e^{Mc_M},\,\,\,\,\,\,\,\,M\rightarrow \infty. \label{06}
\end{equation}
In words, the last term in the above sum cancels all the previous terms. To prove (\ref{06}) observe that
\begin{equation*}
    \sum_{j=1}^{M}e^{jc_j}\leq C e^{\left(M+1\right)c_M},
\end{equation*}
for some positive constant $C$ and
\begin{equation*}
   C e^{\left(M+1\right)c_M}=o\left(e^{\left(M+1\right)c_{M+1}}\right),\,\,\,\,\,\,\,\,M\rightarrow \infty,
\end{equation*}
and the desired result follows.\\

Having information for the quantity $A_N$ is critical regarding our analysis for the asymptotics of relation (\ref{03}).

\textbf{Step II}. On the other hand a substantial amount of information is known for the behaviour of the quantity $L_{1}\left(M;\delta\right)$. In particular, the following dichotomy holds (see, \cite{DPSETS}):
\begin{align*}
&\text{\textbf{(Case I)}}\,\,\,\,\,\,\,\,\,\lim_{M\rightarrow \infty}L_{1}\left(M;\delta\right)<\infty, \,\,\,\,\text{iff}\,\,\,\,\,\sum_{j = 1}^\infty e^{-d_j \xi} < \infty\,\,\,\,\text{for some}\,\, \xi>0,\\
&\text{\textbf{(Case II)}}\,\,\,\,\,\lim_{M\rightarrow \infty}L_{1}\left(M;\delta\right)=\infty, \,\,\,\,\,\text{iff}\,\,\,\,\,\sum_{j = 1}^\infty e^{-d_j \xi} = \infty\,\,\,\,\text{for all}\,\, \xi>0.
\end{align*}
(Case II) has been studied for \textit{rich classes} of sequences, namely those which satisfy conditions ((\ref{cond2})--(\ref{cond3})).
Then, (see, \cite{DPM}, \cite{DPSETS}, and, \cite{DM})
\begin{equation}
  L_{1}\left(M;\delta\right) \sim f(M)\ln\left(\frac{f(M)}{f^{\prime}(M)}\right), \,\,\,\,\,M\rightarrow \infty. \label{08}
\end{equation}
In fact the first five terms of the asymptotic expansion of $L_{1}\left(M;\delta\right)$ are known as $M\rightarrow \infty$. Moreover, these classes were extended in the case where
\begin{equation*}
  d_j=\left(\ln j\right)^{-p},\,\,p>0,
\end{equation*}
see, \cite{CCP LOG} for details. At the end of the paper we present a table including the leading term (asymptotically) for the quantities $\sum_{j=1}^{M}l_j,\,L_{1}(M;\lambda),$ and, $E\left[T_{m}(M;\lambda)\right]$, for several special distributions. \\

\textbf{Conclusion.} Under conditions (\ref{deca})--(\ref{cond3}), and having formula (\ref{03}) (of Theorem 1), as well as powerful tools as the Euler-Maclaurin summation formula, and tangible results for \textit{rich} classes of distributions (see, relation (\ref{08}) and Table 1), the leading asymptotics of $E\left[T_{m}(N;\alpha)\right]$ in the case where the sequence $\alpha$ is the “union” of two subsequences $\beta$ and $\delta$ is easy to compute (as $N\rightarrow \infty$). To obtain these asymptotics one has: (i) to deduce which one of the subsequences $\beta$ and $\delta$ produces the rare set of coupons, (ii) to compute leading asymptotics for the quantity $D_M  +B_M$ (see relation (\ref{03}) and apply the methods described in Step 1 above), and finally, (iii) complete the work by using the asymptotic behaviour of the quantity $L_1\left(M;\delta\right)$ (see, e.g., Table 1 below), where $\delta$ is the subsequence of the rare coupons as we have already indicated. Let us also note, that both distributions contribute to our final answer (see, Subsection 2.1 below).\\\\
\textbf{Remark 2}. Let $r\geq 1$, be an integer. Set
\begin{equation*}
T_{m}\left(N;\alpha\right)^{(r)}:=T_{m}\left(N;\alpha\right)\left(T_{m}\left(N;\alpha\right)+1\right)\left(T_{m}\left(N;\alpha\right)+2\right)\cdots
 \left(T_{m}\left(N;\alpha\right)+r-1\right).
\end{equation*}
Similarly, for the $r$-th rising moment of the random variable $T_{m}\left(N;\alpha\right)$, that is, $E\left[T_{m}\left(N;\alpha\right)^{(r)}\right] $, one has (see, \cite{DPM}, \cite{DPSETS})
\begin{equation*}
E\left[T_{m}\left(N;\alpha\right)^{(r)}\right] =r
 \int_{0}^{\infty}\left[1-\prod_{j=1}^N
 \bigg(1-S_m({p_{j}t})e^{-{p_{j}}t}\bigg)\right]t^{r-1}dt, \quad r=1,2,\cdots,
\end{equation*}
and easily get the leading asymptotics of these moments as $N\rightarrow \infty$. We will give closure to this section with the following (kind of) general
\begin{proposition}
Suppose that we fix the distribution of the half of the coupons.  Then the quantity $E\left[T_{m}\left(N;\alpha\right)^{(r)}\right]$ becomes minimum when the other half of the coupons are uniformly distributed.\end{proposition}
\textit{Proof}. Without lost of generality we may assume that the distribution of the coupons is a mixture of the two distributions generating from the sequences $c_j$ and $b_{j}$. In view of (\ref{bsc}) we have
\begin{align*}
&E\left[T_{m}\left(N;\alpha\right)^{(r)}\right]=r \\
&\times
 \int_{0}^{\infty}\left[1-\prod_{j=1}^M
 \bigg(1-S_m\left({\frac{b_{j}}{A_N}t}\right)\bigg)e^{-{\frac{b_{j}}{A_N}}t}
 \prod_{j=1}^M \bigg(1-S_m\left({\frac{c_{j}}{A_N}t}\right)\bigg)e^{-{\frac{c_{j}}{A_N}}t}
 \right]t^{r-1}dt\\
&=r A_N ^{r} \int_{0}^{\infty}\left[1-\prod_{j=1}^M
 \bigg(1-S_m\left(b_{j}t\right)\bigg)e^{-b_{j}t}
 \prod_{j=1}^M \bigg(1-S_m\left(c_{j}t\right)\bigg)e^{-c_{j}t}
 \right]t^{r-1}dt.
\end{align*}
Let us fix the distribution generated by the sequence $c_{j}$ and set
\begin{equation*}
  h\left(b_1 ,b_2 , \cdots, b_M\right):=\sum_{j=1}^{M}
  \ln\left[\bigg(1-S_m\left(b_j t\right)\bigg)e^{-b_j t}\right].
\end{equation*}
Clearly, the function $h$ is symmetric under all permutations of its variables. Moreover,  it is not hard to see that, for all $i,j=1,2,\cdots,M$, with $i\neq j$:
\begin{equation*}
  \left(b_i
-b_j \right) \left(\frac{\theta h\left(b_1 ,\cdots, b_M\right)}{b_{i}}-\frac{\theta h\left(b_1 ,\cdots, b_M\right)}{b_{j}}\right)\leq 0.
\end{equation*}
Hence, $h$ is a Schur concave function. It follows that $h$ attains its maximum value, when all $b_j$'s (which are strictly positive) are equal, (see, for details \cite{MO}, pp.84 and 413). This completes the proof. \hfill $\blacksquare$\\\\
For a concrete instance of Proposition 2, we refer the reader to Figure 2, where it has been assumed that the total number of coupons ranges from $20$ to $200$. The distribution of the odd coupons is assumed to be fixed; in particular, in view of relation (\ref{sb1}), we take $a_{2j-1}=j$. As for the distribution of the even coupons, it is sequentially assumed that they follow either the uniform distribution ($a_{2j}=1$), the Zipf distribution ($a_{2j}=j^{-2}$), or a positive power law ($a_{2j}=j^{3}$). The number of trials for one complete set ($m=1$) of coupons is smaller when ($a_{2j}=1$).

\subsection{Examples}
Here we illustrate the previous arguments with a few characteristic examples. Notice, that the number $m$ \textit{does not} appear in the leading term of the average $E\left[T_{m}\left(N;\alpha\right)\right]$ as $N\rightarrow \infty$.\\ \\
\textbf{Example 1.} If half of the coupons are coming from the \textit{positive power law} distribution with (positive) parameter $r$ and the rest of the coupons are coming from the \textit{Zipf law},
then,
\begin{equation*}
 \beta:=\{b_j\}_{j=1}^{\infty }=j^{r},\,\,r>0,\,\,\,\,\, \delta:=\{d_j\}_{j=1}^{\infty}=j^{-p},\,\,p>0.
\end{equation*}
In words, the rare coupons are coming from the \textit{Zipf} distribution. Hence, from formulae (\ref{deca}--\ref{cond3}) and the results from Table $1$ we get (recall that $N=2M$)
\begin{equation*}
E\left[T_{m}\left(N;\alpha\right)\right] \sim \frac{N^{r+p+1}}{2^{r+p+1}\left(r+1\right)}\ln N,\,\,\,\,\,\,N \rightarrow \infty.
\end{equation*}
\\
\\
\textbf{Example 2.} Let us consider the case of \textit{exponential - uniform} problem. Then, clearly
\begin{equation*}
 \beta:=\{b_j\}_{j=1}^{\infty }=e^{pj},\,\,p>0,\,\,\,\,\, \delta:=\{d_j\}_{j=1}^{\infty}=1.
\end{equation*}
The sequence $\delta:=\{d_j\}_{j=1}^{\infty}$ does not satisfy conditions (\ref{cond2}--(\ref{cond3})). Hence, we can not apply Theorem 1. Though, using relation (\ref{02}) and the well known results (from Table 1) one immediately has
\begin{equation*}
E\left[T_{m}\left(N;\alpha\right)\right] \sim  \frac{e^{p\left(M+1\right)}}{e^{p}-1}\ln M,\,\,\,\,\,\,M \rightarrow \infty.
\end{equation*}
\\
\textbf{Example 3.} The \textit{uniform versus Zipf distribution} scenario. This problem (for the case of $m=1$) has been studied extensively in (\cite{UVZ}), where the authors derived the first and second moment of $T_{m}\left(N;\alpha\right)$ up to the fifth and sixth (asymptotically) term in order to obtain the corresponding variance as $M\rightarrow \infty$, and the distribution, as well. Following the notation of our paper the rare coupons are coming from the Zipf law. We have
\begin{equation}
E\left[T_{m}\left(N;\alpha\right)\right] \sim
\frac{N^{p+1}\ln N}{2^{p+1} },\,\,\,\,\,\,N \rightarrow \infty, \label{15}
\end{equation}
in accordance with the results of \cite{UVZ}. Figure 1 illustrates the case where half of the coupons are uniformly distributed, while the rest of the coupons follow the standard Zipf distribution. Simulations using Python present $E\left[T_{m}\left(N;\alpha\right)\right]$ in comparison with the theoretical result (see relation (\ref{15})).
\\
\\
\textbf{Example 4.} Suppose now that half of the coupons are coming from the \textit{logarithmic law} with parameter $r>0$ and the rest of the coupons are coming from the \textit{Zipf law}. In this case one has:
\begin{equation*}
 \delta:=\{d_j\}_{j=1}^{\infty }=j^{-p},\,\,p\geq 1,\,\,\,\,\, \beta:=\{b_j\}_{j=2}^{\infty}=\left(\ln j\right)^{-r},\,\,\,\,r>0.
\end{equation*}
Hence, from formula (\ref{03}) and the results from Table $1$ we get (since $N=2M$)
\begin{equation}
E\left[T_{m}\left(N;\alpha\right)\right] \sim  \frac{N^{p+1}}{2^{p+1}\left(\ln N\right)^{r-1}},\,\,\,\,\,\,N \rightarrow \infty.\label{14}
\end{equation}
We will give closure with the following\\

\textbf{Example 5.} Consider the classic coupon collector for the sequences $ \beta:=\{b_j\}_{j=0}^{\infty}=e^{pj}$ and $\delta:=\{d_j\}_{j=0}^{\infty }=e^{-pj}$, where $p>0$, and observe that the elements of these sequences are proportional to each other. Hence, they produce the same coupon probabilities. This follows from the fact that for each $M$, if we let $c_M =e^{-pM}$, then $\left\{b_{j}: 0\leq j \leq M\right\}=\left\{c_M d_{j}: 0\leq j \leq M\right\}$. Thus,
\begin{align*}
E\left[T_{m}\left(M;\beta \right)\right] & = E\left[T_{m}\left(M;\delta \right)\right]\\
& \sim C \,\frac{e^{p\left(M+1\right)}}{e^{p}-1},\quad M\rightarrow \infty,
\end{align*}
where $C$ is a constant. The sequence $\delta$ does not satisfy conditions $(ii)$ and $(iii)$ of relation (\ref{cond3}). Yet by applying relation (\ref{02}) and the well known results (from Table 1) one immediately has (the expected asymptotic result)
\begin{equation}
E\left[T_{m}\left(N;\alpha\right)\right] \sim C \,\frac{e^{p\left(2M+1\right)}}{e^{p}-1},\quad M\rightarrow \infty.\label{14a}
\end{equation}

\section{Generalization}
Let us now consider the case where the family of coupon probabilities is the \textit{interlacing of three subfamilies} one of which consists of far rather coupons (in the sense of Section 1), than the ones of the other two families. For example, and following the notation of Section 1 we assume that $N=3M$ and consider Zipf, exponential, and uniform distributions, namely the sequences $\delta:=\{d_j\}_{j=1}^{\infty }$, $\beta:=\{b_j\}_{j=0}^{\infty }$, and, $\eta:=\{n_j\}_{j=1}^{\infty }$, such that
\begin{align*}
a_{3j}& = b_j=e^{qj},\,\,\,\,q>0,\,\,\,\,\,j=1,2,\cdots,\\
a_{3j+1}&=d_j=j^{-p},\,\,p>1,\,\,\,j=1,2,\cdots,\\
a_{3j-1} &=  n_j=1,\,\,\,\,\,\,\,\,\,\,\,\,\,\,\,\,\,\,\,\,\,\,\,\,\,\,\,j = 1, 2, \cdots.
\end{align*}
Hence, $A_N =B_M +D_M+H_M,$ where
\begin{align*}
D_M:=\sum_{j=1}^M d_j,\,\,B_M:=\sum_{j=0}^{M-1} b_j\,\,\,\,\,\text{and}\,\,\,\,H_M:=\sum_{j=1}^M n_j.
\end{align*}
Clearly, the rare coupons come from the sequence $\delta$, i.e., the Zipf distribution. Using similar arguments as in Section 2 one arrives at the following result
\begin{equation}
E\left[T_{m}\left(N;\alpha\right)\right]=\left(D_M +B_M+H_M\right)\int_{0}^{\infty}\left[1-\prod_{j=1}^M
 \bigg(1-S_m(d_{j}u)e^{-d_{j}u}\bigg)\right]du. \label{20}
\end{equation}
Again, the leading term of the average comes from the integral above, as well as, from the leading term of the quantity: $D_M +B_M+H_M$, which may be evaluated by the Euler-Maclaurin summation formula, and the rest of the tools presented in Section 2. For example, in our case one has
\begin{equation*}
E\left[T_{m}\left(N;\alpha\right)\right]
\sim \frac{e^{q\left(M+1\right)}} {e^{q}-1}\, M^p\ln M,\,\,\, M \rightarrow \infty,
\end{equation*}
where $N=3M$. Observe that only the exponential and the Zipf distributions contribute to the leading term of $E\left[T_{m}\left(N;\alpha\right)\right]$, whereas the number of coupon sets $m$ has no effect on the result.\\
\\
\textbf{Concluding Remarks.} If our family of coupons is the interlacing of $j$ subfamilies, where $j$ is any positive integer, our approach may easily be generalized with similar arguments. Let us also mention that by mimiking this approach one may derive the leading term of the moments of the r.v. (see, relation (\ref{02})). For example, for the second moment of the random variable  $T_{m}\left(N;\alpha\right)$ one has
\begin{align}
  E\left[T_{m}\left(N;\alpha\right)^{2}\right]
  =&  E\left[\sum_{j=1}^{T_{m}\left(M;\delta\right)}W_{j}^{2}\right]+
  2E\left[\sum_{i <j}^{T_{m}\left(M;\delta\right)}W_{i}W_{j}\right] \nonumber\\
  =&E\left[T_{m}\left(M;\delta\right)\right]E\left[W_{j}^{2}\right]+
  E\left[T_{m}\left(M;\delta\right)\left(T_{m}\left(M;\delta\right)-1\right)\right]
  E\left[W_{j}\right]^{2} \nonumber,\label{02b}
\end{align}
and following the same steps one has the (leading) asymptotics as $M\rightarrow \infty$ (see, relation (\ref{geo})). The limiting distribution of $T_{m}\left(N;\alpha\right)$ will be determined mainly by the distribution of the rare coupons. However, in order to obtain the appropriate normalization for $T_{m}\left(N;\alpha\right)$, one must delve deeper into the asymptotics of the first and second moments. In our recent work \cite{UVZ}, where half of the coupons are uniformly distributed and the other half follow the standard Zipf law, and in the case $m=1$,the appropriately normalized random variable $T_{m}\left(N;\alpha\right)$ converges in distribution to a Gumbel random variable. To the best of our knowledge, this is the only published result concerning this variant of the coupon collector problem. We conjecture that, for the classes of distributions studied here, the result under the proper normalization will be the same.\\
Let us also mention that the distribution in the classic version of the problem (the uniform distribution) was first presented by P. Erdős and A. Rényi \cite{E-R}, and, by embedding a Poisson process, by L. Holst \cite{HOLST}. This work and its generalizations to many covering contexts are studied in the paper by Godbole, Grubb, and Kay \cite{G}, where interesting distributional corollaries are provided. Variants of the double Dixie cup problem, like the one presented here, could be applied to several other areas for future research, such as quantum search algorithms, cryptographic engineering, and combinatorial probability (see, e.g., \cite{KESSLER}, \cite{JANA}, and \cite{SAU}).

\section{Table of results and simulation}

Table 1 below presents leading asymptotics for rich classes of distributions and has been used in the examples exhibited in Section 2. Hence, one may compute easily the leading behaviour of the random variable $E\left[T_{m}\left(N;\alpha\right)\right]$ of Theorem 1. Figure 1 exhibits repeated simulations in Python, in the case where half of the coupons follow uniform distribution, while the rest come from the standard Zipf distribution. Simulation versus theoretical results are presented for several values of the total number of coupons $N$. In Figure 2 one inspects the validity of Proposition 2, in the case where the distribution of the odd coupons is fixed. We have assumed that the even coupons follow uniform distribution (Case 1), power law (Case 2), and Zipf law (Case 3). $E\left[T_{m}\left(N;\alpha\right)\right]$ is minimized when the even coupons are uniformly distributed. 

\begin{table}[h!]
\centering
\small
\begin{tabular}{||c c c c c c||}
 \hline
  law& type  & $\sum_{j=1}^{M}l_j$ & $L_{1}(M;\lambda)$ & $E\left[T_{m}(M;\lambda)\right]$ & Ref.\\ [0.5ex]
 $\lambda:=$      &  & (leading  & (leading  & (leading &  \\
 $\left\{l_j\right\}_{j=1}^{\infty}$      &  & term) & term)  & term) &  \\
 \hline\hline
 1 & uniform  &\( M \)  & $\ln M$ & $M \ln M$& \cite{E-R} \\ \\
 $j^{r}$ &  positive & $\frac{M^{r+1}}{r+1}$ & constant  & $C \frac{M^{r+1}}{r+1}$ & \cite{DPSETS}, \cite{DP} \\
  &  power law &  &   &  & \\ \\
 $j^{-p}$& Zipf law &$\frac{M^{1-p}}{1-p}$  & $M^{p}\ln M$  & $\frac{M\ln M}{1-p}$ & \cite{DPSETS}, \cite{DP}, \cite{UVZ}\\
  & $0<p<1$ & &   &  & \\ \\
 $j^{-1}$& std Zipf &$\ln M$  & $M\ln M$  & $M \left(\ln M\right)^{2}$ & \cite{DPSETS}, \cite{DP}, \cite{UVZ}\\
  & $p=1$ & &   &  &  \\ \\

 $j^{-p}$& Zipf law &$\zeta(p)$  & $M^{p}\ln M$  & $\zeta(p)M^{p}\ln M$ & \cite{DPSETS}, \cite{DP}, \cite{UVZ}\\
  & $p>1$ & &   &  &  \\ \\

 $\left(\ln j\right)^{-r}$ & logarithmic & $\frac{M}{\left(\ln M\right)^{r}}$ & $\left(\ln M\right)^{r+1}$  & $M \ln M$ & \cite{CCP LOG}\\ \\

 $e^{pj}$ & exponential & $\frac{e^{p\left(M+1\right)}}{e^{p}-1}$ & constant  & $C\frac{e^{p\left(M+1\right)}}{e^{p}-1}$ &  \cite{DPSETS}, \cite{DP}, \cite{DPM} \\ \\

 $e^{-pj}$  & exponential & $\frac{1}{1-e^{-p}}$ & Example 5  & $C\frac{e^{p\left(M+1\right)}}{e^{p}-1}$ &  \cite{DPSETS}, \cite{DP}, \cite{DPM}\\ \\

 $e^{ic_j},$   & super &    &  &  & \\
  $c_{j}	\nearrow \infty$  & exponential & $e^{M c_{M}}$ & constant  & $Ce^{M c_{M}}$ & \cite{DPSETS}, \cite{DM}, \cite{DPM} \\ \\
 $j!$ & factorial & $M!$ & constant  & $Ce M!$ & \cite{DPSETS}, \cite{DP}, \cite{DPM}\\ \\

 $\frac{1}{j!}$ & reciprocal  & $e$ & $C_{1}M!$  & $C M!$ & \cite{DPSETS}, \cite{DP}, \cite{DPM}\\

   & factorial  &   & &  & \\ \\ [1ex]

 \hline
\end{tabular}
\caption{Leading asymptotics for $\sum_{j=1}^{M}l_j,\,L_{1}(M;\lambda),$ and, $E\left[T_{m}(M;\lambda)\right]$. The parameters $p$ and $r$ are strictly positive. As usual, we use $C$ and $C_{1}$ to denote constants in any equation.}
\label{table:1}
\end{table}
\normalsize

\begin{figure}[H]
\vspace{-22mm}
\centering
\includegraphics[scale =0.5]
{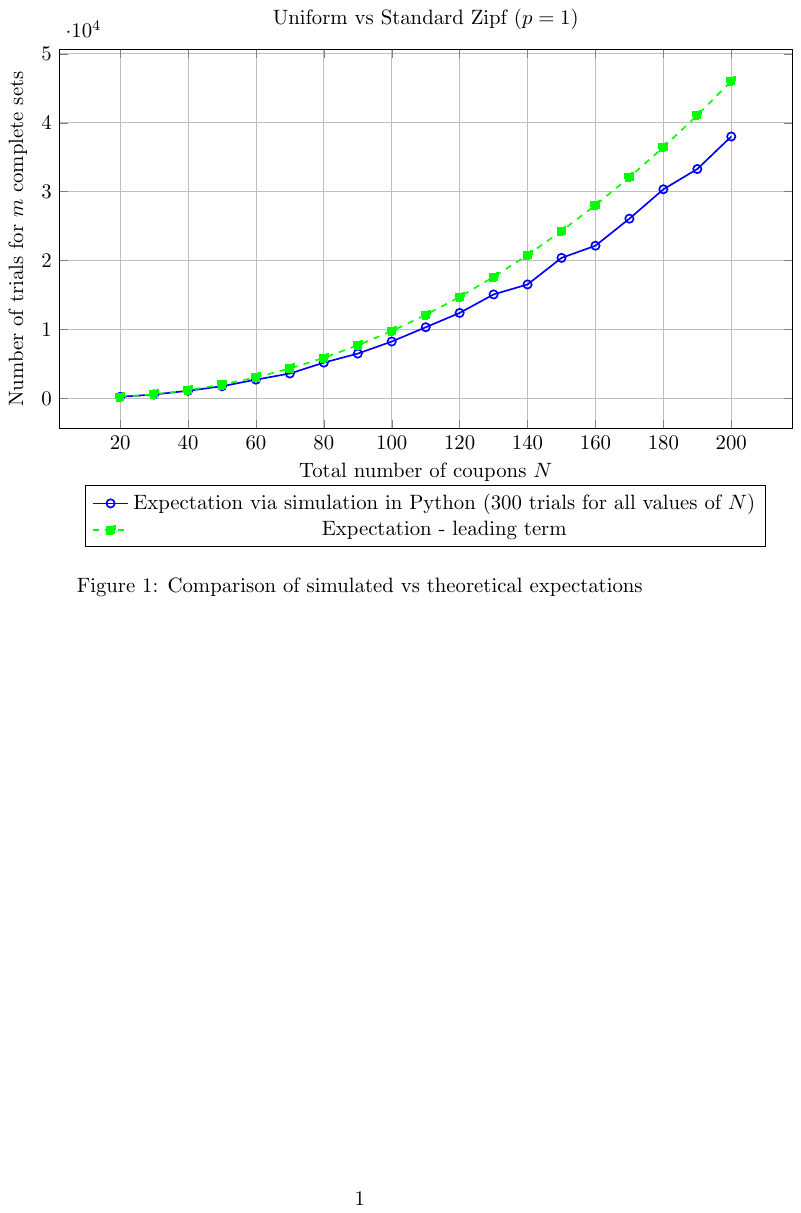
}
\end{figure}

\textbf{Acknowledgements.} 
The author has been partially supported by project MIS 5154714 of the National Recovery and Resilience Plan Greece 2.0 funded by the European Union under the NextGenerationEU Program.


\end{document}